\documentclass{article}

\usepackage[english]{babel}
\usepackage{amsmath, amssymb, amsthm, amsfonts, bm}
\usepackage{enumerate}
\usepackage{graphicx, color}
\usepackage[colorlinks=true]{hyperref}
\usepackage{enumitem}
\usepackage{csquotes}
\numberwithin{equation}{section}
\newtheorem{theorem}{Theorem}[section]
\theoremstyle{definition}

\theoremstyle{remark}

\newcommand{\Rnum}[1]{\uppercase\expandafter{\romannumeral #1\relax}}

\newcommand{\mb}[1]{\mathbb{#1}}
\newcommand{\mc}[1]{\mathcal{#1}}

\DeclareMathOperator*{\essinf}{essinf}

\def\clap#1{\hbox to 0pt{\hss#1\hss}}

\title{A note on the Mukenhoupt class $A_1({\mb{R}})$}
\author{Eleftherios N. Nikolidakis, Andreas G. Tolias}
\date{\today}

\begin{document}
\maketitle

\begin{abstract}
We provide a direct proof of the best possible reverse H\"{o}lder inequality 
 %%%%%%%%%%%%%%%%%%%%%%%%%%   which is satisfied for 
 satisfied by  
  %%%%%%%%%%%%%%%%%%%
 every weight defined on the interval $({0,1})$ with $A_1$-constant equal to $c>1$.
\end{abstract}

\section{Introduction} \label{sec:0}
 
An integrable function $\varphi\colon({0,1})\rightarrow \mb{R}^{+}$ is called an $A_1$-weight if there exists $c>1$ such that:
\begin{equation} \label{eq:1p1}
\frac{1}{|{I}|}\int_{I}\varphi(x)dx\leq c\cdot\essinf_{x\in I}\varphi({x})\,,  
\end{equation} 
for every subinterval $I$ of $({0,1})$. 

The 
  minimum 
 %%%%%   least
  constant $c$ for which \eqref{eq:1p1} holds for all such intervals $I$ is called the $A_1$-constant of $\varphi$ and it is
   denoted 
%%%%%%   defined 
by $[{\varphi}]_{A_1}$.

An equivalent definition of an $A_1$-weight is given as follows: 
For every $\varphi \in L^{1}\big({({0,1})}\big)$ the maximal function $\mc{M}\varphi$ is defined by 
\begin{equation} \label{eq:1p2}
	\mc{M}\varphi(x) = \sup\left\{ \frac{1}{|I|} \int_I |\varphi|\; :\;  x\in I,\; 
	I\;\text{is an open subinterval of}\;({0,1}) \right\},
\end{equation}
for every  $x\in ({0,1})$. It is easily seen that  \eqref{eq:1p1} is equivalent to the   condition: 
\begin{equation}
 \label{eq:1p3}
\mc{M}\varphi(x) \leq c\,\varphi(x)\,,\;\text{for almost every}\; x\in({0,1})\,.
\end{equation}
We denote by $A_1\big({(0,1),c}\big)$ the class of all positive functions on $({0,1})$ which satisfy \eqref{eq:1p1} or
 equivalently satisfy \eqref{eq:1p3}. 

In a similar manner one can define the class $A_1\big({J,c}\big)$ for any interval $J\subset\mb{R}$ and any $c>1$, as the collection of all functions $\varphi\colon J\rightarrow \mb{R}^{+}$ which satisfy \eqref{eq:1p1} for all subintervals $I$ of $J$. Respectively for  $\varphi \in L^{1}({J})$, $\mc{M}_{J}\varphi\colon J\rightarrow \mb{R}^{+}$ is defined 
by the rule \eqref{eq:1p2}, where the supremum is now taken over all open subintervals of $J$ that contain the point $x$. Then $\varphi\in A_1\big({J,c}\big)$ if and only if $\mc{M}_{J}\varphi\leq c\,\varphi$, almost everywhere on $J$. 

 For $\varphi$, $J$ as above,  the quantity $[{\varphi}]_{A_1}^{J}$,  which is the $A_1$ constant of $\varphi$ on the interval $J$, is defined analogously to the case $J=(0,1)$. 
  Consequently $[{\varphi}]_{A_1}$ denotes $[{\varphi}]_{A_1}^{({0,1})}$ and $\mc{M}\varphi$ stands for $\mc{M}_{({0,1})}\varphi$, for any nonnegative $\varphi \in L^{1}\big({(0,1)}\big)$. 

A notion that is connected to the class $A_1\big({(0,1),c}\big)$ is that of the decreasing rearrangement of $\varphi$, which is denoted by $\varphi^{*}$, where $\varphi^{*}\colon({0,1})\rightarrow \mb{R}$ is the unique decreasing, right continuous  function on $(0,1)$, equimeasurable to $\varphi$. The function $\varphi^{*}$ can also be given by a specific formula (see for example   \cite{1}, page 39.)

It is  proved in \cite{2}, that for any $\varphi \in A_1\big(((0,1),c\big)$, the function $\varphi^{*}$ also belongs to $A_1\big((0,1),c\big)$, that is it satisfies the condition
\begin{equation}\label{eq:1p4}
\frac{1}{t} \int_0^t \varphi^*(u) \, du \leq c \, \varphi^*(t), \quad\mbox{ for every } t \in (0,1).
\end{equation}

Moreover, in \cite{2}   the following is proved

\begin{theorem}
If $\varphi\colon (0,1) \rightarrow \mathbb{R}^+$ is a function belonging to the class $A_1\big((0,1),c\big)$, then the reverse H\"{o}lder inequality
\begin{equation}\label{eq:1p5}
\frac{1}{|I|} \int_I \varphi^p(x) dx\leq \frac{1}{c^{p-1}(c+p - c\,p)} \left({ \frac{1}{|I|} \int_I \varphi(x)dx }\right)^p,	
\end{equation}
is satisfied for every interval $I \subset (0,1)$ and for every 
$p \in \big[{1, \frac{c}{c-1}}\big)$. 

Moreover the constant in the right hand of \eqref{eq:1p5} is the best possible such that the inequality  
    is satisfied for every  weights $\varphi$ with $A_1$-constant equal to $c$    and all $p$ with $1\le p<\frac{c}{c-1}$.

%Moreover, inequality \eqref{eq:1p5} is the best possible for any such $p$, if one considers all weights $\varphi$ with %$A_1$-constant equal to $c$.\qed
\end{theorem}

The above theorem in fact implies that the best possible range of values of $p$, for which an $A_1$-weight with $A_1$-constant $[{\varphi}]_{A_1}$ equal to $c$ is $L^p$-integrable on $(0,1)$, is the interval $\big[{1,\frac{c}{c-1}}\big)$.

For the above two results, alternative proofs have been given in \cite{4}. The approaches for proving the above theorem  given in \cite{2} and \cite{4}, use the notion of the decreasing rearrangement of a weight $\varphi$.

In this short note we provide a direct  proof of the above theorem without passing through the  notion of 
decreasing rearrangement. 
 Lastly, we should mention that the dyadic version of the above problem has been studied and completely solved in \cite{3}.

\section{Proof of the Theorem}
Let $\varphi\colon(0,1)\rightarrow\mathbb{R}^+$ be an integrable function such that $\mc{M}\varphi \leq c\,\varphi$, almost everywhere on $(0,1)$.  
We let $f := \int_0^1 \varphi(x)dx$ and for every  we define $\lambda > 0$ the set  
\[
E_\lambda = [\mc{M}\varphi > \lambda] = \big\{ {x\in(0,1)\;:\; \mc{M}\varphi(x) > \lambda }\big\}.\]  

Note that for $\lambda \in ({0,f})$ we have that $E_\lambda=(0,1)$.  
Suppose now that $\lambda$ satisfies $\lambda \geq f$. By definition we have that the set $E_\lambda$ is open, thus it can be written as an at most countable disjoint union of open subintervals of $(0,1)$. Thus  we may write
$E_\lambda = \bigcup\limits_{j\in D_\lambda} I_j$, where the set of indices $D_\lambda$ is at most countable, each $I_j$ is an open subinterval of  $(0,1)$ and the family $(I_j)_{j\in D_{\lambda}}$ is pairwise disjoint. We consider the following two cases:

\begin{enumerate}
\item[(i)] 
$(I_j)_{j\in D_\lambda} = \{({0,1})\}$, so that $E_{\lambda}= (0,1)$. In this case we have that 
\[
|{E_\lambda}| = 1 \geq \frac{f}{\lambda} =  \frac{\int_{E_\lambda} \varphi(x)dx}{\lambda}\,.
\]
\item[(ii)] 
The family $(I_j)_{j\in D_\lambda}$ differs from $\{(0,1)\}$. Then for each $j\in D_{\lambda}$ and each $\delta > 1$ close enough to $1$, one can define an open subinterval $I_{j,\delta}$ of $(0,1)$ for which    $I_j \subset I_{j,\delta}$ and
 $|I_{j,\delta}| = \delta |I_j|$.

Then, for each such $j$ and $\delta$,  $I_{j,\delta}$ is not contained in  $E_\lambda$,  
 thus there exists a point $x_{j,\delta} \in I_{j,\delta}$ for which $\mc{M}\varphi(x_{j,\delta}) \leq \lambda$. Hence, by definition of the maximal function $\mc{M}\varphi$, we get that:
\[\frac{1}{|I_{j,\delta}|} \int_{I_{j,\delta}} \varphi(x)dx\leq \lambda, 
\]
  for every $\delta > 1$ sufficiently close to  $1$. Letting $\delta \to 1^{+}$ in the above inequality it follows
that \[
\frac{1}{|I_j|} \int_{I_j} \varphi (x)dx\leq \lambda.
\]
Since this holds for every $j\in D_{\lambda}$ and 
 $ (I_j)_{j\in D_{\lambda}}$ are pairwise disjoint   with their  union equal to $E_\lambda$ we obtain that
 $\frac{1}{|E_\lambda|} \int_{E_\lambda} \varphi(x)dx \leq \lambda$
 which implies that
\begin{equation}\label{eq:2p1}
|E_\lambda| \geq \frac{1}{\lambda}\int_{E_\lambda} \varphi (x)dx\,.
\end{equation}
\end{enumerate}
We have just proved that  inequality \eqref{eq:2p1} holds for any $\lambda \geq f = \int_0^1 \varphi(x)dx$.

Consider   now any $p$ in the interval $\big[{1, \frac{c}{c-1} }\big)$.
 Multiplying both parts of \eqref{eq:2p1} with $p \lambda^{p-1}$ and integrating on the interval $[f, +\infty)$ we obtain:
\begin{equation}\label{eq:2p2}
\int_{f}^{+\infty} p\,\lambda^{p-1} \big| E_\lambda\big| \,d\lambda \geq 
\int_{ f}^{+\infty} p\,\lambda^{p-2} \Big( \int_{E_\lambda}\varphi(x)dx \Big)\,d\lambda\,. 
\end{equation}
By  Fubini's theorem and taking into account that $E_\lambda=[\mc{M}\varphi>\lambda]$    we get that 
\begin{equation}\label{eq:2p3}
\int_0^1 \big(\mc{M}\varphi (x)\big)^pd x =\int_{0}^{+\infty} p\,\lambda^{p-1} \big| E_\lambda \big| \,d\lambda.
\end{equation}

Using  \eqref{eq:2p3} and \eqref{eq:2p2} and taking into account that $E_\lambda=(0,1)$ when $0<\lambda< f$ we get that 
\begin{eqnarray}
\int_0^1 \big(\mc{M}\varphi (x)\big)^pdx &= & \int_{0}^{f}  p\,\lambda^{p-1} \big|   E_\lambda\big|d\lambda +\int_{f}^{+\infty} p\,\lambda^{p-1} \big|   E_\lambda\big| \,d\lambda\nonumber    \\
    & \geq & \int_{0}^{f}p\,\lambda^{p-1}\,d\lambda +\int_{f}^{+\infty} p\,\lambda^{p-2} \Big( {\int_{E_\lambda} \varphi(x)dx} \Big)\,d\lambda\nonumber    \\
  &= & f^p +\int_{f}^{+\infty} {\int_{E_\lambda}  p\,\lambda^{p-2}    \varphi}(x)dx  \,d\lambda\,.\label{eq:2p4}
\end{eqnarray}
We apply now Fubini's Theorem to   the integral on the right-hand side of \eqref{eq:2p4}.
\begin{align}
& \int_{f}^{+\infty} \int_{E_\lambda}  p\,\lambda^{p-2}    \varphi(x)dx  \,d\lambda 
 =\int_{E_f}  \varphi(x) \int_{f}^{\mc{M}\varphi(x)}p\,\lambda^{p-2}d\lambda\, dx \nonumber\\
 & =\int_{E_f}  \varphi(x)  \frac{p}{p-1} \Big( (\mc{M}\varphi(x))^{p-1}-f^{p-1}\Big) \,d\lambda \, dx\nonumber\\
 &= \int_{0}^{1} \frac{p}{p-1}   \varphi(x) \Big( \big(\mc{M}\varphi(x)\big)^{p-1}-f^{p-1}\Big)d\lambda \, dx\nonumber\\
 &=-\frac{p}{p-1} f^p  +  \frac{p}{p-1} \int_0^1 \varphi(x) \big(\mc{M}\varphi(x)\big)^{p-1}dx \label{eq:2p5}
\end{align}
where we took into account that $\mc{M}\phi(x)\ge f $ for all $x\in (0,1)$.

From \eqref{eq:2p4} and \eqref{eq:2p5} we get that
\[\int_0^1 \big(\mc{M} \varphi(x)\big)^pdx  \geq -\frac{f^p}{p-1} + 
 \frac{p}{p-1}\int_0^1 \varphi(x) (\mc{M}\varphi(x))^{p-1} \]
 thus
 \begin{equation}
 \int_0^1 \Big\{ p \, \varphi(x)\,(\mc{M}\varphi(x))^{p-1} - (p-1)\,(\mc{M} \varphi(x))^p \Big\} dx 
\leq f^p \label{eq:2p6}
 \end{equation}

Consider now for each    $y > 0$, the function
$h_y\colon [{y, +\infty})\rightarrow\mathbb{R}$ given by the formula:
\[
h_y(t) = p \,y \, t^{p-1} - (p-1)\,t^p.
\]
We  notice that
\[
h_y'(t) = p(p-1)\, y \, t^{p-2} - p(p-1)\, t^{p-1} = p(p-1)\, t^{p-2} (y - t) < 0 
\]
for all $t>y$,
thus $h_y$ is strictly decreasing on the interval $[{y, +\infty})$.
Hence for every $t$ with 
 $y \leq t \leq c\, y$, we   have that
\begin{equation}
h_y(t) \geq h_y(c\, y) = y^p\,c^{p-1}  (c+p-cp).   \label{ineqh}
\end{equation} 

We also notice  that 
\[  \varphi(x) \leq \mathcal{M}\varphi(x) \leq c\, \varphi(x)  \mbox{ for almost every } x\in (0,1)\]
where the left hand inequality follows by Lebesgue's differentiation theorem for $L^1(0,1)$-functions,
while the right hand inequality follows from our assumption that  the function $\varphi$ belongs to $A_1\big((0,1),c\big)$.
 Applying \eqref{ineqh}   for $y=\varphi(x)$ and $t=\mc{M}\varphi(x)$, we get that for almost all $x\in (0,1)$
 it holds that 
\begin{equation}
p \, \varphi(x)\,(\mathcal{M}\varphi(x))^{p-1} - (p-1)(\mathcal{M} \varphi(x))^p \geq
  \big(\varphi(x)\big)^p\, c^{p-1}\left( c+p-cp \right) . \,\label{eq:2p7}
\end{equation}

Combining now \eqref{eq:2p6} and  \eqref{eq:2p7}, we deduce that:
\begin{align}
c^{p-1}\left(c+p-cp\right) \int_0^1 \varphi^p(x) dx
\leq f^p=\Big(\int_0^1 \varphi(x)dx\Big)^p   \,.\label{eq:2p8}
\end{align}
Since $p \in \big[{1, \frac{c}{c-1}}\big)$, we have that  $c+p-cp > 0$,  therefore from \eqref{eq:2p8} we get that
\[
\int_0^1 \varphi^p(x)dx \leq 
\frac{1}{c^{p-1} ({c + p - c\,p})} \Big({ \int_0^1 \varphi(x)dx }\Big)^p\,,
\]
which is the inequality that is stated in \eqref{eq:1p5}, for  
the whole interval $I = (0,1)$. Now if one considers an arbitrary interval $I \subset (0,1)$, then, by defining $\varphi_I := \varphi | I \colon I \to \mathbb{R}^+$, we easily  
see that we   have that
\[
\mathcal{M}_I \varphi_I \leq c \,\varphi_I, \quad \text{a.e. in } I.
\]
Using then the same   arguments as in the above proof,  
but now working on the interval \( I \) instead, we get inequality  
\eqref{eq:1p5} for this specific \( I \), concluding the proof of the  
reverse H\"{o}lder inequality.

Finally, we prove the sharpness of \eqref{eq:1p5} as follows:  
Let $c > 1$ and  consider the function $\varphi \colon (0,1) \to \mathbb{R}^+$ which is given by
\[
\varphi(t) = t^{-1 + \frac{1}{c}}, \quad t \in (0,1)\,.
\]
Then $\int_0^1 \varphi(t)dt= c$, and $\frac{1}{t} \int_0^t \varphi(u)\, du = c \,\varphi(t)$, for all $t \in (0,1)$, and since $\varphi$ is decreasing on $(0,1)$, we get that $\varphi$ satisfies: $[\varphi]_{A_1} = c$\,

Moreover,  $\int_0^1 \varphi^p(t)dt = \frac{1}{1 + ( -1+\frac{1}{c})p}$, \; for any   $p \in \big[1, \frac{c}{c-1} \big)$. Thus for any given $p$ in the aforemention interval we   have that:
\[
\frac{\int_0^1 \varphi^p(t)dt }{\big({ \int_0^1 \varphi(t)dt }\big)^p}
= \frac{1}{c^p \left(1 +  ( -1+\frac{1}{c})p \right) }
 = \frac{1}{c^{p-1}(c + p - c\,p)}
\]
which gives us the sharpness  
of \eqref{eq:1p5}  for any $p \in \big[1, \frac{c}{c-1} \big)$\,.\qed

\vspace{20pt}
\noindent Nikolidakis Eleftherios, Associate  Professor, Department of Mathematics, University of Ioannina, 45110, Greece\\
E-mail address: enikolid@uoi.gr

\vspace{15pt}
\noindent Tolias Andreas, Associate Professor, Department of Mathematics, University of Ioannina, 45110, Greece \\
E-mail address: atolias@uoi.gr

\end{document}